\documentclass{article}
\usepackage{amsmath}
\usepackage{amssymb}

\usepackage{parskip}

\newcommand{\abs}[1]{\left|{#1}\right|}

\usepackage{ifpdf}
\ifpdf
\usepackage[hypertex,hyperindex,pdftex,breaklinks=true]{hyperref}%
\hypersetup{%
pdfauthor={Joerg Arndt},%
pdftitle={On computing the generalized Lambert series},%
pdfsubject={Lambert series},%
pdfkeywords={Lambert series, Theta convergence, fast computation, inverse Fibonacci numbers},%
bookmarksopen,
bookmarksopenlevel={0},
}
\else
\usepackage[hypertex,hyperindex,breaklinks=true]{hyperref}%
\fi

\begin{document}
\bibliographystyle{plain}
\title{On computing the generalized Lambert series}
\author{J\"{o}rg Arndt}
\date{
Georg Simon Ohm University of Applied Sciences\\
 Nuremberg, Germany
\\[5mm]
\texttt{<arndt@jjj.de>}
\\[5mm]
\today
}
\pagestyle{plain} \maketitle

\begin{abstract}\noindent
We show how the generalized Lambert series
$\sum_{n\geq{}1}{{x q^n}/{(1-x q^n})}$
can be computed with Theta convergence.
This allows the computation of the sum of the
inverse Fibonacci numbers without splitting the sum into even and odd part.
The method is a special case of an expression
for the more general series $\sum_{n\geq{}0}{{t^n}/{(1-x q^n)}}$,
which can be obtained from either the Rogers-Fine identity
or an identity by Osler and Hassen.

Addendum (June 2012):
the main identity \eqref{rel:lambert-qxt} already appeared
in Agarwal's 1993 paper \cite{agarwal-lambert-ram}.
This is now remarked at the end of this paper.
\end{abstract}

\section{The identity}

Fine's function $F(a,b,t)$ can be given \cite[p.1]{Fine} as
\[
F(a,b,t) =
 \sum_{n\geq{}0}{
 \frac
    {(1-a\,q)\,(1-a\,q^2)\,(1-a\,q^3)\,\cdots\,(1-a\,q^{n})}
    {(1-b\,q)\,(1-b\,q^2)\,(1-b\,q^3)\,\cdots\,(1-b\,q^{n})}
 \, t^n }
\]
Defining $(a;q)_0=1$ and $(a;q)_n=(1-a)\,(1-a\,q)\,(1-a\,q^2)\,\cdots\,(1-a\,q^{n-1})$
(for $n\in\mathbb{N}_+$) we can write
\[
F(a,b;t) =
 \sum_{n\geq{}0}{
 \frac{(a\,q;q)_n}{(b\,q;q)_n} \, t^n }
\]
The Rogers-Fine identity \cite[relation (14.1), p.15]{Fine}
is
\[
(1-t)\,F(a,b;t) =
\sum_{n\geq{}0}{ \frac{(a\,q;q)_n\,(a\,t\,q/b;q)_n}{(b\,q;q)_n\,(t\,q;q)_n}
 \, (1-a\,t\,q^{2n+1})\,b^n\,t^n\,q^{n^2} }
\]
Replacing $b$ by $a\,q$ gives, for the left side,
\[
(1-t)\,\sum_{n\geq{}0}{ \frac{(a\,q;q)_n}{(a\,q^2;q)_n} \, t^n } =
(1-t)\,\sum_{n\geq{}0}{ \frac{1-a\,q}{1-a\,q^{n+1}} \, t^n }
\]
For the right side we get
\[
\sum_{n\geq{}0}{ \frac{(a\,q;q)_n\,(t;q)_n}{(a\,q^2;q)_n\,(t\,q;q)_n}
 \, (1-a\,t\,q^{2n+1})\,(aq)^n\,t^n\,q^{n^2} } =
\]
\[
\sum_{n\geq{}0}{ \frac{(1-a\,q)\,(1-t)}{(1-a\,q^{n+1})\,(1-t\,q^n)}
 \, (1-a\,t\,q^{2n+1})\,(aq)^n\,t^n\,q^{n^2} }
\]
Now we divide both sides by $(1-a\,q)\,(1-t)$ to obtain
\[
\sum_{n\geq{}0}{ \frac{t^n}{1-a\,q^{n+1}} } =
\sum_{n\geq{}0}{ \frac{1}{(1-a\,q^{n+1})\,(1-t\,q^n)}
 \, (1-a\,t\,q^{2n+1})\,(aq)^n\,t^n\,q^{n^2} }
\]
Finally replacing $a$ by $x/q$ gives the desired identity
\begin{equation}\label{rel:lambert-qxt}
\sum_{n\geq{}0}{ \frac{t^n}{1-x\,q^{n}} } =
\sum_{n\geq{}0}{
 \frac{ (1-x\,t\,q^{2n}) }{(1-x\,q^{n})\,(1-t\,q^n)} \, x^n\, t^n\, q^{n^2} }
\end{equation}
Both sides are convergent for $\abs{q}<1$, $\abs{x}<1$, and $\abs{t}<1$.
The right side reveals symmetry in $x$ and $t$ which is also observed
as relation (4.3) in \cite[p.655]{Osler}:
\begin{equation}\label{rel:symm}
\sum_{n\geq{}0}{ \frac{t^n}{1-x\,q^{n}} } =
\sum_{n\geq{}0}{ \frac{x^n}{1-t\,q^{n}} }
\end{equation}

Setting $x=q$, $t=q$ in \eqref{rel:lambert-qxt}, and multiplying by $q$ gives
the Clausen's \cite[p.95]{Clausen}
\begin{equation}\label{rel:clausen-lambert}
L(q) =
\sum_{n\geq{}1}{ \frac{q^n}{1-q^{n}} } =
\sum_{n\geq{}1}{ \frac{1+q^{n}}{1-q^{n}} \, q^{n^2} }
\end{equation}
We denote the rate of convergence of the series on the right side
as \emph{Theta convergence}. (This is faster than linear but not what
is usually called super-linear convergence,
which should really be called exponential convergence.)

\section{Summation of inverse Fibonacci numbers}
The computation of $\sum_{n\geq{}1}{{1}/{F_n}} \approx 3.359885\ldots$ where
$F_0=0$, $F_1=1$, and $F_n=F_{n-1}+F_{n-2}$ for $n\geq{}2$
(the Fibonacci sequence) has been considered in various publications
(see, for example,
 \cite{Bruckman},
 \cite{Castellanos},
 \cite{Gould},
 and \cite{Horadam}).
The following identities for the sums with even and odd indices
have been known for a long time.
%
Let $\alpha=1/2\,(1+\sqrt{5})$, $\beta=-1/\alpha=1/2\,(1-\sqrt{5})$,
then (even indices, see \cite[p.468]{Knopp})
\[
\sum_{n\geq{}1}{ \frac{1}{F_{2n}} } =
\sqrt{5}\,\left( L(\beta^2) - L(\beta^4) \right)
 \quad\approx 1.5353705\ldots
\]
and (odd indices, see \cite[pp.94]{PiAGM})
\[
\sum_{n\geq{}1}{\frac{1}{F_{2n-1}}} =
 \frac{\sqrt{5}}{4} \, \left( \Theta_3^2(\beta^2) - \Theta_3^2(\beta)\right)
 \quad\approx 1.8245151\ldots
\]
where $\Theta_3(q)=1+2\,\sum_{n\geq{}1}{q^{n^2}}$.
Both relations allow computation with Theta convergence.
Other expressions for the sums are
\[
\sum_{n\geq{}1}{ \frac{1}{F_{2n}} } =
\sum_{n\geq{}1}{ \frac{\beta^{\,n(n+1)}}{ 1 - \beta^{2n} } }
\]
and
\[
\sum_{n\geq{}1}{\frac{1}{F_{2n-1}}} =
 -\sqrt{5} \,  \beta \, \left( \sum_{n\geq{}0}{ \beta^{\,2n(n+1)}}  \right)^2
\]

The search for a relation that does not require the splitting between even and
odd indices (and works for more than just the Fibonacci numbers) has been the
main motivation for this work.
Gosper gave the following formula for
the Fibonacci sequence \cite[p.66]{AIM304}:
\begin{equation}\label{rel:gosper-stunt}
\sum_{n\geq{}1}{ \frac{1}{F_{n}} } =
 \sum_{n\geq{}0}{
 \frac{(-1)^{n(n-1)/2} \, \left(F_{4n+3}+(-1)^n\,F_{2n+2}\right)}
      {F_{2n+1}\,F_{2n+2}\; G_1\,G_3\,G_5\,\cdots\,G_{2n+1} } }
\end{equation}
where $G_n=2\,F_{n-1} + F_n$ (Lucas numbers).
Note the last term in the denominator
 has been corrected from $G_{2n-1}$ to $G_{2n+1}$.

We will need a specialization of \eqref{rel:lambert-qxt},
often called \emph{generalized} Lambert series.
First set $t=q$ and multiply by $x$ to obtain
\[
\sum_{n\geq{}0}{ \frac{x\,q^n}{1-x\,q^{n}} } =
\sum_{n\geq{}0}{
 \frac{ (1-x\,q^{2n+1}) }{(1-x\,q^{n})\,(1-q^{n+1})} \, x^n\, q^{n^2+n} }
\]
Replace $x$ by $x\,q$ to obtain the desired form
(note the new limits of summation)
\begin{equation}\label{rel:lambert-gen}
L(x,q) =
\sum_{n\geq{}1}{ \frac{x\,q^n}{1-x\,q^{n}} } =
\sum_{n\geq{}1}{
 \frac{ (1-x\,q^{2n}) }{(1-x\,q^{n})\,(1-q^{n})} \, x^n\, q^{n^2} }
\end{equation}
%
%

The computation of $\sum_{n\geq{}1}{{1}/{f_n}}$ where
$f_0\geq{}0$, $f_1\geq{}1$, and $f_n=m_1\,f_{n-1}+m_2\,f_{n-2}$ for $n\geq{}2$
has been treated by Horadam \cite{Horadam}.
We restrict our attention to cases where $f_0=0$, $f_1=1$ (Fibonacci type sequences),
$m_1\geq{}1$, $m_2\neq{}0$, and $\Delta=m_1^2+4\,m_2 > 0$.
The special case $m_1=m_2=1$ corresponds to the Fibonacci numbers.
%
The roots of the characteristic polynomial $x^2-(m_1\,x+m_2)$
are $\alpha$ and $\beta$ where
\[
 \alpha=1/2\,\left(m_1+\sqrt{\Delta}\right), \qquad
 \beta=1/2\,\left(m_1-\sqrt{\Delta}\right)
\]
The following identity, given by Horadam
 \cite[relation (4.17), p.108]{Horadam},
together with relation \eqref{rel:lambert-gen}, is the key to
fast computation of $\sum_{n\geq{}1}{{1}/{f_n}}$:
\begin{equation}\label{rel:fast-invfib}
\sum_{n\geq{}1}{\frac{1}{f_n}} =
 (\alpha - \beta)\, \left( \frac{1}{\alpha-1} + L\left( \frac{1}{\alpha},\, \frac{\beta}{\alpha}\right) \right)
\end{equation}
%
%
%
%
%
%

\section{A few more relations}

We mention several identities observed during this work.

\subsection{Alternate relation for fast computation}
Fine's second transformation is \cite[relation (12.2), p.13]{Fine}
\[
 (1-t)\,\sum_{n\geq{}0}{
 \frac{(a\,q;q)_n}{(b\,q;q)_n} \, t^n } =
%
\sum_{n\geq{}0}{ \frac{(b/q;q)_n}{(b\,q;q)_n\,(t\,q;q)_n}
 \, (-a\,t)^{n}\,q^{(n^2+n)/2} }
\]
Setting $b=a\,q$, then replacing $a$ by $a/q$ and dividing by $(1-a)\,(1-t)$,
and finally replacing $a$ by $x$ gives an alternate identity for fast computation:
\begin{equation}\label{rel:lambert-qxt-alt}
\sum_{n\geq{}0}{ \frac{t^n}{1-x\,q^{n}} } =
%
 \sum_{n\geq{}0}{
 \frac{ (q;q)_n }{ (x;q)_{n+1} \, (t;q)_{n+1} } \, (x\,t)^n\, q^{(n^2-n)/2} }
\end{equation}
The symmetry between $x$ and $t$ is evident from the right side,
giving \eqref{rel:symm} once more.

\subsection{Another derivation of the main identity}
Identity \eqref{rel:lambert-qxt} can alternatively be deduced from
the (remarkable!) relations (6.5), (6.6), and (6.7)
by Osler and Hassen \cite[p.659]{Osler}, namely
\begin{subequations}
\begin{eqnarray}
\label{rel:osler-12}
\sum_{n\geq{}0}{ \frac{ \alpha^n \, q^{d\,(a\,n+b)} }{ 1 - \beta\,q^{c\,(a\,n+b)}} } =
\sum_{n\geq{}0}{ \frac{ \beta^n \, q^{b\,(c\,n+d)} }{ 1 - \alpha\,q^{a\,(c\,n+d)} } } =
%
\\
\label{rel:osler-3}
 = \sum_{n\geq{}0}{ \left(
 \frac{ \alpha^n \, \beta^n \, q^{(a\,n+b)\,(c\,n+d)} }{ 1 - \beta\,q^{c\,(a\,n+b)} } +
 \frac{ \alpha^{n+1} \, \beta^n \, q^{(a\,(n+1)+b)\,(c\,n+d)} }{ 1 - \alpha\,q^{a\,(c\,n+d)} }
 \right) }
\end{eqnarray}
\end{subequations}
by setting $a=c=1$, $b=d=0$, and replacing $\alpha$ by $x$ and $\beta$ by $t$.
%

\subsection{The symmetric relation \eqref{rel:symm} from an identity by Fine}
%
%
Fine's relation (16.3) \cite[p.18]{Fine} is
\[
F(a,b,t) =
\frac{ (a\,q;q)_{\infty} }{ (b\,q;q)_{\infty} } \,
  \sum_{n\geq{}0}{ \frac{(b/a;q)_n}{(q;q)_n} \, \frac{(aq)^n}{ 1 - t\,q^n } }
\]
The following can be obtained by setting $b=a\,q$, dividing by $(1-a\,q)$,
and finally replacing $a$ by $x$:
\[
\sum_{n\geq{}0}{ \frac{t^n}{ 1 - x\,q^{n+1} } } =
\sum_{n\geq{}0}{ \frac{(xq)^n}{ 1 - t\,q^n} }
\]
Replacing $x$ by $x/q$ we get relation \eqref{rel:symm} once again.

\subsection{More symmetric relations}
An identity resembling \eqref{rel:symm} is
\[
\sum_{n\geq{}0}{ \frac{t^n}{(x;q)_{n+1}} } =
\sum_{n\geq{}0}{ \frac{x^n}{(t;q)_{n+1}} }
\]
This follows from an identity
given by Gosper \cite{Gosper-mathfun}
\[
\frac{ (t;q)_\infty \, (x;q)_\infty }{ (q;q)_\infty} \,
 \sum_{n\geq0}{ \frac{t^n}{(x;q)_{n+1} } }
=
\sum_{n\geq{}0}{ \frac{ (t;q)_n \, (x;q)_n }{ (q;q)_n } \, q^n }
\]
and noting that both the right side and the factor in front of the sum
on the left side are symmetric in $x$ and $t$.

Setting $A=B=C=D=1$ in relations \eqref{rel:osler-12} and \eqref{rel:osler-3}
(and again replacing $\alpha$ by $x$ and $\beta$ by $t$) we obtain
\[
 t \, \sum_{n\geq{}1}{ \frac{ x^n\,q^n }{ 1 - t \, q^n } }
 \; = \;
 x \, \sum_{n\geq{}1}{ \frac{ t^n\,q^n }{ 1 - x \, q^n } } =
\]
\[
 = \sum_{n\geq{}1}{ \frac{ x^n\,t^n\,q^{n^2} }{ 1 - t\, q^n } } +
   x \, \sum_{n\geq{}1}{ \frac{ x^n\,t^n\,q^{n(n+1)}}{ 1 - x \, q^n } }
\]
\[
 = \sum_{n\geq{}1}{ \frac{ x^n\,t^n\,q^{n^2} }{ 1 - x\, q^n } } +
   t \, \sum_{n\geq{}1}{ \frac{ x^n\,t^n\,q^{n(n+1)}}{ 1 - t \, q^n } }
\]
%
\[
 = \sum_{n\geq{}1}{
  \frac{ (1-x\,t\,q^{2n}) }{(1-x\,q^{n})\,(1-t\,q^n)} \, x^n\, t^n\, q^{n^2} }
\]
The last form reveals the symmetry in $x$ and $t$.

\subsection{Another derivation of the relation for the generalized Lambert series}
Relation \eqref{rel:lambert-gen} can be obtained from the following relation
 given by Knuth (attributed to J.\ R.\ Wrench, Jr.)
 \cite[p.644, solution to exercise 5.2.3-27]{DEK3}:
\[
\sum_{n\geq{}1}{ \frac{a_n\,q^n}{ 1 - q^n } } =
\sum_{n\geq{}1}{ \left[ a_n + \sum_{k\geq{}1}{ (a_n+a_{n+k})\,q^{k\,n} } \right]\,q^{n^2} }
\]
In the following (well-known) relation
\[
\sum_{n\geq{}0}{ \frac{x\,q^n}{ 1 - x\,q^n } } =
\sum_{n\geq{}1}{ \frac{x^n}{ 1 - q^n } }
\]
replace $x$ by $q\,x$ to get
\[
L(x,q) =
\sum_{n\geq{}1}{ \frac{x\,q^n}{ 1 - x\,q^n } } =
\sum_{n\geq{}1}{ \frac{x^n\,q^n}{ 1 - q^n } }
\]
Now set $a_n=x^n$ to obtain
\[
L(x,q) =
\sum_{n\geq{}1}{ \left[ 1 + \sum_{k\geq{}1}{ (1+x^k)\,q^{k\,n} } \right]\,x^n\,q^{n^2} } =
 \sum_{n\geq{}1}{ \left[ \sum_{k\geq{}0}{ q^{k\,n} } + \sum_{k\geq{}1}{ x^k\,q^{k\,n} } \right]\,x^n\,q^{n^2} }
\]
which gives identity \eqref{rel:lambert-gen} in the form
\[
L(x,q) =
 \sum_{n\geq{}1}{ \left[ \frac{1}{1-q^n} + \frac{x\,q^n}{1-x\,q^n}\right] \, x^n \, q^{n^2} }
\]
%

\subsection{Gosper's matrix-product for Clausen's identity}
%
Gosper defines \cite{Gosper-mathfun} matrices $K(k,n)$ and $N(k,n)$ as
\[
K(k,n)=
\begin{bmatrix}
q^{n+2k+1} & {q\,(1-q^{2k+n})}/{\left((1-q^k) \, (1-q^{k+n})\right)} \\
 0 & 1
\end{bmatrix}
\]
\[
N(k,m)=
\begin{bmatrix}
q^{k} & {q}/{(1-q^{k+n})} \\
0 & 1
\end{bmatrix}
\]
These matrices satisfy
\[
 N(k,n) \cdot K(k,n+1) = K(k,n) \cdot N(k+1,n)
\]
Writing
 $K(k,\infty)$ for $\lim_{n\to\infty}{K(k,n)}$ and
 $N(\infty,n)$ for $\lim_{k\to\infty}{N(k,n)}$
we obtain \eqref{rel:clausen-lambert} as upper right entries on both sides of
\[
\prod_{n\geq{}0}{N(1,n)} \cdot \prod_{k\geq{}1}{K(k,\infty)} =
\prod_{k\geq{}1}{K(k,0)} \cdot \prod_{n\geq{}0}{N(\infty,n)}
%
\]

\subsection{Jordan's identity for the bilateral series}

For the bilateral series (\emph{Jordan-Kronecker Function}),
defined as
\[
f(x,t) =
\sum_{n=-\infty}^{+\infty}{ \frac{t^n}{ 1 - x\, q^n } }
\]
we have $f(x,t)=f(t,x)$ and
the following two relations given in \cite[chap.3]{VC-ellram}
\[
f(x,t) =
%
\sum_{n=-\infty}^{+\infty}{ q^{n^2} \,x^n\,t^n\, \left( +1 + \frac{x\,q^n}{1-x\,q^n} + \frac{t\,q^n}{1-t\,q^n} \right) }
\]
\[
f(x,t) =
\sum_{n=-\infty}^{+\infty}{ q^{n^2} \,x^n\,t^n\, \left( - 1 + \frac{1}{1-x\,q^n} + \frac{1}{1-t\,q^n} \right) }
\]
As pointed out by Shaun Cooper,
the latter is already given in \cite[rel.(22), p.450]{Jordan2}.
Both forms are equivalent to
\begin{equation}\label{rel:jordan}
f(x,t) =
\sum_{n=-\infty}^{+\infty}{ \frac{(1-t\,x\,q^{2n})}{ (1-x\,q^n)\,(1-t\,q^n)} \, x^n\,t^n\,q^{n^2} }
\end{equation}
%

\section*{Addendum: 1993 appearence of the main identity}
As of June 2012 it has come to the author's attention that
the main identity \eqref{rel:lambert-qxt} already appeared
in Agarwal's 1993 paper \cite{agarwal-lambert-ram}.
Agarwal proves identity \eqref{rel:jordan},
which is his relation (6.16),
and gives our identity \eqref{rel:lambert-qxt} as his relation (6.17).

%


\end{document}